\documentclass[a4paper,11pt]{article}
\usepackage[latin1]{inputenc}
\usepackage{amsmath, amssymb, amsthm, amsfonts, bbm}
\usepackage{pifont,fontenc,tabularx}
\usepackage[english]{babel}
\usepackage{graphicx,xcolor}
\usepackage{multirow}
\usepackage[compact]{titlesec}

\theoremstyle{definition}

\usepackage{tikz}

\allowdisplaybreaks

\newcommand{\da}{\,\mathrm{d}a}
\newcommand{\dxi}{\,\mathrm{d}\xi}

\newcommand{\M}{\mathcal{M}}
\newcommand{\R}{\mathbb{R}}

\newcommand{\LBO}{\Delta_{\M}}
\newcommand{\surfGrad}{\nabla_{\M}}
\newcommand{\diver}{\operatorname{div}}

\newcommand{\RR}{\mathbb R}
\newcommand{\ZZ}{\mathbb Z}
\newcommand{\II}{\mathcal I}
\newcommand{\LL}{L}

\newcommand{\spa}{\operatorname{span}}

\newcommand{\cell}{\triangle}

\renewcommand{\S}{\mathbf{S}}

\newcommand{\I}{\mathcal{I}}

\definecolor{light-gray1}{gray}{0.60}
\definecolor{light-gray2}{gray}{0.75}
\definecolor{light-gray3}{gray}{0.90}

\title{On Isogeometric Subdivision Methods for PDEs on Surfaces}
\author{Bert J\"uttler,Angelos Mantzaflaris,Ricardo Perl,Martin Rumpf}

\begin{document}




\maketitle

{ \bf Abstract:} 
Subdivision surfaces are proven to be a powerful tool in geometric modeling and
computer graphics, due to the great flexibility they offer in capturing irregular topologies. 
This paper discusses the robust and efficient implementation of 
an isogeometric discretization approach to partial differential equations on surfaces using subdivision methodology.
Elliptic equations with the Laplace-Beltrami and the surface bi-Laplacian operator as well as the associated eigenvalue problems are considered. 
Thereby, efficiency relies on the proper choice of a numerical quadrature scheme which preserves the expected higher order consistency.
A particular emphasis is on the robustness of the approach in the vicinity of extraordinary vertices. 
In this paper, the focus is on Loop's subdivision scheme on triangular meshes.
Based on a series of numerical experiments, different quadrature schemes are compared and a mid-edge quadrature, which is
easy-to-implement via lookup tables, turns out to be a preferable choice due to its robustness and efficiency. 

\vspace{2ex}

\section{Introduction} \label{sec1}


During the last years, \emph{isogeometric analysis} (IgA)
\cite{HuCoBa05,cottrell:2009} emerged as a means of unification of the
previously disjoint technologies of geometric design and numerical
simulation. The former technology is classically based on parametric
surface representations, while the latter discipline relies on finite
element approximations. The advantages of the isogeometric paradigm 
include the elimination of both the geometry approximation
error and the computation cost of changing the representation between
design and analysis. 
A particular challenge is the problem of efficient numerical integration.
Indeed, the higher degree and multi-element
support of the basis functions seriously affect the cost of robust numerical integration.
Consequently, numerical quadrature in IgA is an active area of research
\cite{Antolin2015817, HuReSa10, MaJ15final, MaJueKh14, Schillinger20141}.

Subdivision schemes are widespread in geometry processing and computer graphics. 
For a comprehensive introduction to subdivision methods in general we refer the reader to \cite{PeRe08} and \cite{Ca12}. 
With respect to the use of subdivision methods in animation see \cite{DeKaTr98, ThWaSt06}.
Nowadays, subdivision finite elements are also extensively used 
in engineering \cite{CiOrSc00, CiLo11, CiScAn02, CiOr01, GrTuSt02, GrTu04}.
This way, the subdivision methodology changed the modeling paradigms in \emph{computer-aided design} (CAD)
systems over the last years tremendously, e.g. the \emph{Freestyle} module of PTC Creo\circledR, the \emph{PowerSurfacing} add-on for SolidWorks\circledR~or \emph{NX RealizeShape} of Siemens PLM\circledR.
For the integration of subdivision methods in CAD systems, we refer to 
\cite{AnBeCa13, BoZo04} and the references therein.  
With this development, powerful \emph{computer-aided engineering} (CAE)
codes have to be implemented for subdivision surfaces. Today NURBS
and subdivision surfaces co-exist in hybrid systems.

Among the most popular subdivision schemes are the Catmull-Clark
\cite{CaCl78} and Doo-Sabin \cite{DoSa78} schemes on quadrilateral meshes, and
Loop's scheme on triangular meshes \cite{Lo94}.
The Loop subdivision basis functions have been extensively studied, most 
interestingly regarding their smoothness \cite{Re95, ReSc01}, (local) linear independence \cite{PeWu06a, ZoJKo14}, approximation power \cite{Ar01}, and the
robust evaluation around so called extraordinary vertices \cite{St99}.

Computing PDEs on surfaces is an indispensable tool either to approximate physical processes on the surfaces \cite{GrKrSc02} 
or to process textures on the surface or the surface itself \cite{BaXu03, ClDiRu04, HiPo04}.
Recently, surface PDEs have been considered in the field of IgA, and
both numerical and theoretical advancements are reported \cite{DeQu15, LaMaMo14, LaMo14}.

The use of Loop subdivision surfaces as a finite element discretization technique is discussed in \cite{CiOrSc00} and isogeometric discretizations
based on the Catmull-Clark scheme in \cite{Ba13a, WaHiPo11}, whereas \cite{DiSnOr12}
investigates the use of the Doo-Sabin scheme in the finite element context.
An isogeometric finite element analysis based on Catmull-Clark
subdivision solids can be found in \cite{BuHaUm10}.

\begin{figure}[h]
\begin{center}
\includegraphics[scale=1.]{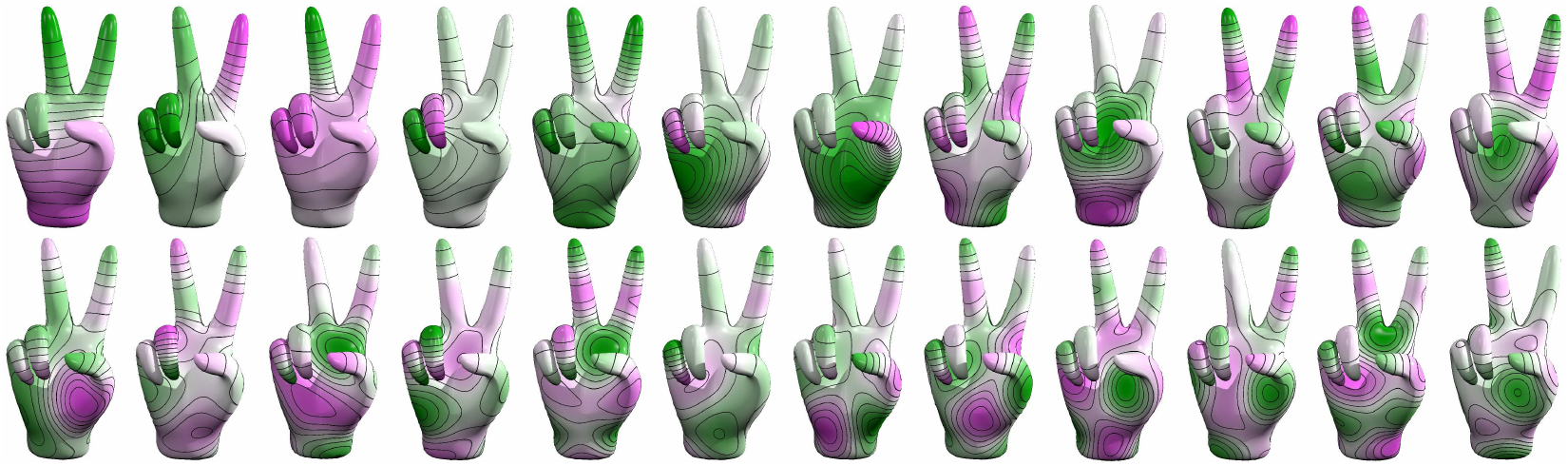}
\end{center}
\caption{As an application we show here the first 24 eigenfunctions of the Laplace-Beltrami operator computed on a complex surface using the isogeometric subdivision approach. 
See Figure \ref{fig:handB} for the coarse control mesh.}
\label{fig:handB_eigen}
\end{figure}

In this paper, we focus on the Loop subdivision scheme \cite{Lo87} on triangular meshes. 
Loop subdivision is known to describe 
$C^2$ limit surfaces $\M$ except at finitely many points (called extraordinary vertices)
where they are only $C^1\cap H^2$ \cite{Re95, ReSc01}.
This allows us to use them in a conforming finite element approach not only for second-order but also for fourth-order elliptic PDEs.
As model problem we consider the \emph{Laplace-Beltrami} equation
\begin{align}
 - \LBO u &= f \text{  on } \M\,, \label{eq::strongLaplaceProblem}
\end{align}
the \emph{surface bi-Laplacian} equation
\begin{align}
 (- \LBO)^2 u &= f \text{  on } \M \,  \label{eq::strongBiLaplacianProblem}
\end{align}
as well as the eigenvalue problem
\begin{align}
 -\LBO u = \lambda u  \text{  on } \M \,.  \label{eq::strongEigenvalProblem}
\end{align}
Let us remark that on closed, smooth surfaces the kernel of the Laplace-Beltrami and the surface bi-Laplacian operators are the constant functions.  

The crucial step in the implementation of a subdivision finite element method is the choice of numerical quadrature. 
The \emph{computational cost}, \emph{consistency}, \emph{robustness} and the \emph{observed order of convergence} are important
parameters to evaluate the appropriateness of each numerical integration technique.
We aim at a detailed experimental study of the convergence behavior for different quadrature rules.
We examine Gaussian quadrature, barycenter
quadrature, mid-edge quadrature and discuss adaptive strategies around extraordinary
vertices. Furthermore, we provide a look-up table for the mid-edge quadrature that facilitates the implementation 
and results in a fairly robust simulation tool and a very efficient assembly of finite element matrices. 
This allows for a straightforward extension of existing subdivision modeling codes to simulations with subdivision surfaces for applications, see Figure \ref{fig:handB_eigen}.

\emph{Organisation of the paper.}
In Section \ref{sec2}, we fix some notation and define subdivision functions and subdivision surfaces. 
The isogeometric subdivision method to solve elliptic partial differential equations on surfaces is described in Section \ref{sec3}, where we  
derive the required finite element mass and stiffness matrices.
Section \ref{sec4} investigates the different numerical quadrature rules and comments on the efficient implementation. 
In Section \ref{sec5}, we report on the observed convergence behavior for a set of selected test cases.
Finally, we draw conclusions in Section \ref{sec6}.

\section{Subdivision Functions and Surfaces}  \label{sec2}

To discretize geometric differential operators and to solve geometric partial differential equations on closed subdivision surfaces 
we consider {\em isogeometric function spaces} defined via {\em Loop's subdivision}.

The domain of a Loop subdivision surface is a topological manifold
obtained by gluing together copies of a standard triangle. More
precisely, we consider the unit triangle $\triangle\subset\RR^2$
with vertices $(0,0)$, $(1,0)$ and $(0,1)$, which is
considered as a closed subset of $\RR^2$, and a finite {\em cell index
  set} $\II_c\subset\mathbb\ZZ$. The Cartesian product
$\triangle\times\II_c$ defines the {\em pre-manifold} which consists
of the {\em cells} $\triangle\times\{i\}$, $i\in\II_c$.  The manifold
is constructed by identifying points on the boundaries of the cells,
as described below.

In addition to the cell index set, we consider an {\em edge index set}
$\II_e\subset\mathbb\II_c\times\II_c$, which is assumed to be {\em
  symmetric}
\begin{equation*}
(i,j)\in\II_e\Rightarrow(j,i)\in\II_e
\end{equation*}
and {\em irreflexive} 
\begin{equation*}
\forall i\in\ZZ: \; (i,i)\not\in\II_e.
\end{equation*}
Each edge is described by a symmetric pair of edge indices $((i,j)$,
$(j,i))\in\II_e\times\II_e$.  Moreover, each cell contributes to
exactly three edges,
\begin{equation*}
\forall i\in\II_c: |\{j: (i,j)\in\II_e\}|=3. 
\end{equation*}

\begin{figure}[h]
\begin{center}
\includegraphics[scale=1.]{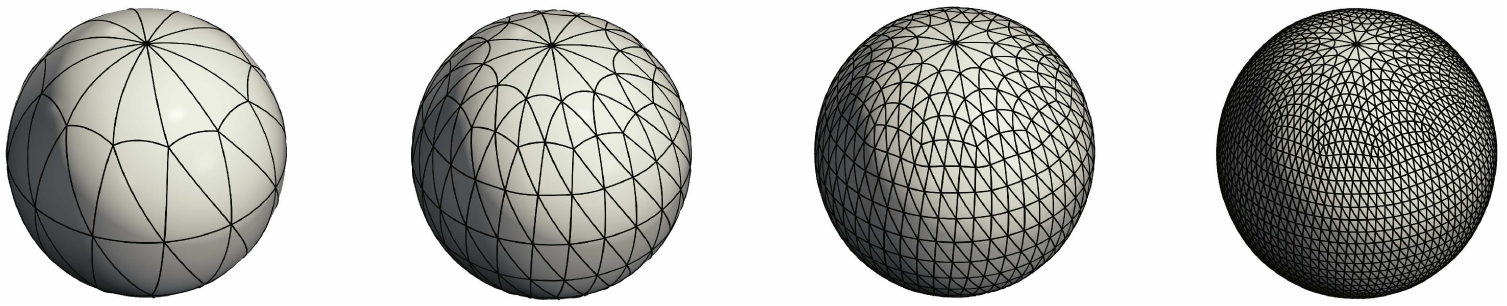}
\end{center}
\caption{The spherical subdivision limit surface Spherical-5-12 with successively refined control meshes.
The control meshes are plotted on the limit surface, i.e., they visualize the curved IgA-elements on the surface.
Furthermore, the control meshes possess two vertices of valence $12$
and $24$ vertices of valence $5$ on every refinement level.
Since two EVs share an edge in the coarsest control mesh, this
mesh cannot be used for simulation.  }
\label{fig:SphereControl}
\end{figure}

Each edge index $(i,j)$ is associated with a displacement
$\delta^i_j$, which is one of the 18 isometries ($3\cdot 3 \cdot 2$) that map the unit
triangle to one of its three neighbor triangles obtained by reflection and index permutation. 
These displacements in particular identify the common edges of neighboring triangles and its orientation
with respect to the two triangles. They need to satisfy the following
two conditions.  Firstly, the identification of the common edge and
its orientation has to be consistent for each edge, thus
\begin{equation*}
\forall (i,j)\in\II_e:\; \delta^i_j=(\delta^j_i)^{-1}. 
\end{equation*}
Secondly, each edges of a triangle is identified with exactly one
edge of another triangle,
\begin{equation*}
\forall (i,j)\in\II_e:\; \forall(i,k)\in\II_e: \;
    \delta^i_j(\Delta)\cap\Delta=\delta^i_k(\Delta)\cap\Delta
\Rightarrow j=k.
\end{equation*}

Two points of the pre-manifold are identified, denoted by
$(\xi,i)\sim(\eta,j)$, if $(i,j)$ is an edge index and the
displacement $\delta^i_j$ transforms $\eta$ into $\xi$. This implies
that $\xi$ and $\eta$ are located on the boundary of the standard
triangle. We denote with $\hat\sim$ the reflexive and transitive
closure of this relation. This closure leads to the obvious identification of common vertices. 
The topological manifold
$(\triangle\times\II_c)/\hat\sim$
is the {\em domain manifold} of the Loop subdivision surface.

There, the indexed standard triangles define an initial triangulation (of level $0$)
of the domain manifold. Now, Loop's subdivision scheme iteratively creates triangulations of level $\ell>0$
containing $4^\ell|\II_c|$ triangles. They are obtained by creating new
vertices at the edge midpoints and splitting each triangle of level
$\ell-1$ into four smaller ones. The valence of the vertices in these
triangulations are always $6$, except for those vertices of
the coarsest level that possess a different valence. The latter ones
are called {\em extraordinary vertices} (EVs). 
For an example of a Loop subdivision surface with extraordinary vertices and different resolutions of the control mesh see Figure \ref{fig:SphereControl}.

Next, we consider the spaces $\LL^\ell$ of piecewise linear functions on the
triangulation of level $\ell$. Each function is uniquely described by
its {\em nodal values}, i.e., by its values at the vertices of the
triangulation. The {\em Loop subdivision operators}
$R^\ell:\LL^{\ell-1}\to\LL^\ell$ are linear operators that transform a
piecewise linear function $f^{\ell-1}$ of level $\ell-1$ into a
function of level $\ell$. The nodal values of the latter (refined)
function are computed as weighted averages of nodal values of the
coarser function at neighboring vertices.  Values at newly created
vertices (at edge midpoints) depend on the four nodal values of the
coarser function on the two triangles intersecting at the edge, 
while values at existing vertices are computed from
the values at the (valence+1)-many vertices in the one-ring
neighborhood. The weights are determined by simple formulas \cite{SDS96}.

The space of {\em Loop subdivision splines} \cite{Lo87} consists of the limit functions
generated by the subdivision operators,
\begin{equation*}
\{\lim_{L\to\infty}R^LR^{L-1}\cdots R^1f : f\in \LL^0\}. 
\end{equation*}
These limit functions are known to be bivariate quartic polynomials on
all triangles of the triangulation of level $\ell$ that do not possess
any EVs. They are $C^2$-smooth everywhere, except at the EVs.
If one uses a special parameterization around EVs, which is given by
the charateristic map \cite{Re95, ReSc01} then the limit
functions are $C^1$-smooth at EVs.

Let $\II_v$ be the {\em vertex index set} of the coarsest triangulation. Each
vertex $i \in\II_v$ has an associated piecewise linear hat function
$\Lambda_i\in\LL^0$ which takes the nodal value 1 at the associated
vertex and 0 else. The limit functions
\begin{equation*}
\Phi_i=\lim_{L\to\infty}R^LR^{L-1}\cdots R^1\Lambda_i,
\end{equation*}
which will be denoted as the {\em Loop basis functions}, span the
space of Loop subdivision splines. The support of the Loop basis
function is the two-ring neighborhood of the vertex $i$ in the
coarsest triangulation. If the support does not contain EVs in its
interior, then it consists of $24$ triangles of the coarsest level and
the Loop basis function is the $C^2$-smooth quartic box spline on the
type-1 triangulation. Otherwise, the Loop basis function is still a piecewise polynomial
function but possesses an infinite number of polynomial segments. 
These are associated with triangles obtained by successive refinement 
in the vicinity of the EV and not containing the EV itself \cite{PeRe08}.

Now we have everything at hand to introduce a {\em Loop subdivision surface}
\begin{align*}
X:\;  \cell\times\I_c/\hat\sim \to \R^3: \; 
X( \xi , k ) = \sum\limits_{i \in \I_v} C_i \varPhi_i(\xi,k)\,,
\end{align*}
which is obtained by assigning control points $C_i\in\RR^3$ to the
Loop basis functions with $\xi = (\xi_1, \xi_2) \in \cell$. This surface is $C^2$ smooth everywhere, except
at EVs (in the same sense as described before)
\cite{Re95, ReSc01}. Hence, it can be shown that the surface possesses a well-defined
tangent plane at EVs, provided that the control points are not in
a singular configuration. We denote the surface (as an embedded
manifold) by $\M = \M[X]$.
In the context of isogeometric analysis, the parameterization $X$ of the
domain is referred to as the {\em geometry mapping}. The basis
functions $\varphi_i$ of the associated {\em isogeometric function
  space} are the push-forwards of the subdivision basis functions
\begin{align*}
\varphi_i(x) = \varPhi_i \circ X^{-1}(x),
\end{align*}
where $x = X(\xi,k)$, i.e., the basis functions are defined by the discretization $X$. 
We define the discretization space as
\begin{align*}
 V_h = \Big\{ u_h \in \spa\limits_{i \in \I_v} \{ \varphi_i \} \Big\}. 
\end{align*}
Here, $h$ indicates the grid size of the control mesh, i.e., $h= \max_{(i,j) \in \II_e} |C_i-C_j|$.
For each cell, we denote by 
$J( \xi , k ) = \begin{pmatrix} X_{,1} & X_{,2} \end{pmatrix}$ the Jacobian of $X$ and by $G = J^TJ$ its first fundamental form,
where $X_{,1}$ and $X_{,2}$ denote the tangent vectors with $X_{,l} = \frac{\partial }{\partial \xi_l}X( \xi , k )$.
Furthermore, for a function $u: \M \to \R$
\begin{align*}
\surfGrad u(x) = \left(J G^{-T} \nabla_\xi (u \circ X)\right)( \xi , k )
\end{align*}
is the (embedded) tangential gradient or surface gradient and 
\begin{align*}
\LBO u(x) = \frac{1}{\sqrt{\det G}} \diver_\xi\left( \sqrt{\det G} G^{-1} \nabla_\xi (u \circ X) \right) ( \xi , k )
\end{align*}
the Laplace-Beltrami operator, where $x = X( \xi , k )$.

\section{Isogeometric Subdivision Method to solve PDEs}  \label{sec3}

Similar to the classical \emph{finite element method}, the \emph{isogeometric subdivision approach} 
transforms the strong formulation of a partial differential equation 
(e.g. \eqref{eq::strongLaplaceProblem} --
\eqref{eq::strongEigenvalProblem})
into a corresponding weak formulation on a suitable subspace $V^0$ of some Sobolev space 
and approximates the solution of the weak problem in the finite-dimensional sub-spaces $V^0_h \subset V^0$.
We consider 
$$
V^0=H^{1}(\M) \cap \{ \int_{\M} v \da = 0\}
$$
for problem \eqref{eq::strongLaplaceProblem} 
and 
$$
V^0=H^{2}(\M) \cap \{ \int_{\M} v \da = 0\}
$$
for problem  \eqref{eq::strongBiLaplacianProblem}.
Following the general isogeometric paradigm we use  
$$V^0_h = V_h \cap \{ \int_{\M} v_h \da = 0\}
$$
as the discrete ansatz space on the subdivision surface $\M$, e.g., see Figure \ref{fig:SphereControl}.

Let $a: V^0 \times V^0 \to \R$ denote a symmetric, coercive and bounded bilinear form on $V^0$ and $\ell_f: V^0 \to \R$ a linear 
form for a given function $f \in L^2(\M)$. Then there exists a unique solution $u\in V^0$ of the variational problem 
\begin{align}
a(u,v) = \ell_f(v), \quad \forall~v \in V^0, \label{eq:generelWeakForm}
\end{align}
where $\ell_f(v) = \int_{\M} f v \da$ (\cite{BrSc02,Br07,Ci02,Dz88}).
For our model problems, we use
$$
a(u,v) = \int_{\M} \nabla_{\M} u \cdot \nabla_{\M} v \da
$$
and 
$$
a(u,v) = \int_{\M} \Delta_{\M} u \,\Delta_{\M} v \da
$$
for the Laplace-Beltrami \eqref{eq::strongLaplaceProblem} 
and the Bi-Laplacian problem \eqref{eq::strongBiLaplacianProblem}, respectively.
The associated  Galerkin approximation asks for the unique solution $u_h \in V^0_h$ of the discrete variational problem
\begin{align}
 a(u_h,v_h) = \ell_f(v_h), \quad \forall~ v_h \in V^0_h\,. \label{eq:generelDiscreteWeakForm}
\end{align}
Due to the known $H^{2}$-regularity  of the Loop subdivision splines this Galerkin approximation is conforming, i.e., $V^0_h \subset V^0$.
Using the basis expansion $u_h = \sum\limits_{i\in\I_v} u_i \varphi_i$ for $u_h$ with coefficients $u_i\in \R$ 
one obtains the linear system 
\begin{align}
 \S \mathbf{U} = \mathrm{B}, \label{eq:exactLinSystem}
\end{align}
where $\mathbf{U} = (u_i)_{i \in \I_v}\in \R^{\mid\I_v\mid}$ denotes the coefficient vector, $\S_{ij} = a(\varphi_i,\varphi_j) \in \R^{\mid\I_v\mid\times\mid\I_v\mid}$
the stiffness matrix, and $\mathrm{B}_j = \int_\M f \varphi_j \da \in \R^{\mid\I_v\mid}$ the right-hand side.
The stiffness matrix for the Laplace-Beltrami problem \eqref{eq::strongLaplaceProblem} is given by
\begin{align}
\S_{ij}^{\Delta}  = \int_\M \nabla_\M \varphi_i \cdot \nabla_\M \varphi_j \da = \sum\limits_{k \in \I_c} \int_{\cell} \nabla_\xi \varPhi_i G^{-T} \nabla_\xi \varPhi_j \sqrt{\det G} \dxi\label{eq:bilinearLaplace}
\end{align}
and for the surface bi-Laplacian problem \eqref{eq::strongBiLaplacianProblem}  we obtain 
\begin{align}
 \S_{ij}^{\Delta^{\!2}} &=  \sum\limits_{k \in \I_c} \int_{\cell} \diver_\xi\left( \sqrt{\det G} G^{-1} \nabla_\xi \varPhi_i \right) 
 \diver_\xi\left( \sqrt{\det G} G^{-1} \nabla_\xi \varPhi_j \right) \frac1{\sqrt{\det G}} \dxi\,. \label{eq:bilinearBiLaplacian}
\end{align}
The variational formulation of the eigenvalue problem \eqref{eq::strongEigenvalProblem} consists in finding a solution $(u,\lambda) \in V^0 \times \R$ such that
\begin{align*}
 a(u,v) = \lambda \cdot m(u,v), \quad ~ \forall ~ v \in V^0,
\end{align*}
where $m(u,v) = \int_{\M} u \, v \da$ denotes the $L^2$-product on $\M$. Furthermore, we denote by $\mathbf{M}$ the mass matrix with 
\begin{align}
\mathbf{M}_{ij} = m(\varphi_i,\varphi_j) = \sum\limits_{k \in \I_c} \int_{\cell} \varPhi_i \varPhi_j \sqrt{\det G} \dxi. \label{eq:massMatrix}
\end{align}
We obtain the discrete eigenvalue problem $\S \mathbf{U} = \lambda_h \mathbf{M} \mathbf{U}$ 
which can be solved by \emph{inverse vector iteration with projection} \cite{ScWe92a}.

\section{Numerical Quadrature}  \label{sec4}
In this section, we recap the evaluation-based assembly of the previously defined discrete variational formulations.
We review standard Gaussian and barycentric quadrature assembly, discuss the special treatment at  
extraordinary vertices and finally introduce an edge-centered assembly strategy which is based on simple fixed,  
geometry-independent lookup tables and already very good consistency in our experiments (see Section \ref{sec5}).
As discussed above, the isogeometric subdivision approach is based on higher order spline discretizations.
The integrands of the corresponding IgA-matrices (\eqref{eq:bilinearLaplace} -- \eqref{eq:massMatrix}) are in general nonlinear 
and even on planar facets of the subdivision surface (with constant metric $G$) high-degree polynomials. Indeed, apart from EVs the limit functions 
$\Phi_i$ are quartic polynomials and thus for the mass matrix the integrand is a polynomial of degree $8$ and 
for the stiffness matrix \eqref{eq:bilinearLaplace} of degree $6$ and for \eqref{eq:bilinearBiLaplacian} of degree $4$.

In the variational formulation, the integration is always pulled-back to the reference triangle, i.e., 
$$
\int_\M f \da = \sum\limits_{k\in\I_c} \int_\triangle f \circ X \, \sqrt{\det G}(\xi,k) \dxi\, .
$$ 
Thus numerical quadrature is applied to these pulled-back integrands.
For a general function $g$ using the evaluation-based quadrature 
$$
\int_\triangle g(\xi) \mathrm{d} \xi \approx \sum\limits_{q=1}^{K} w_{q} g(\xi^q) \, ,
$$
with quadrature points $\xi_q$ on the reference triangle 
$\triangle$ and weights $w_{q}$, we replace the stiffness matrices, the mass matrix, and the right-hand side 
by the following quadrature based counterparts
\begin{align}
\mathbf{\tilde S}_{ij}^{\Delta} &=& \sum\limits_{k \in \I_c} \sum\limits_{q=1}^{K} w_{q} \left(\nabla_\xi \varPhi_i^T G^{-T} \nabla_\xi \varPhi_j \sqrt{\det G}\right)(\xi^q,k) \label{eq:modifiedBilinearLaplace}
\end{align}
and
\begin{align}
\mathbf{\tilde S}_{ij}^{\Delta^{\!2}} &=& \sum\limits_{k \in \I_c} \sum\limits_{q=1}^{K} w_{q} \left( \diver_\xi\left( \sqrt{\det G} G^{-1} \nabla_\xi \varPhi_i \right) \diver_\xi \left( \sqrt{\det G} G^{-1} \nabla_\xi \varPhi_j \right)  \frac1{\sqrt{\det G}} \right)(\xi^q,k) \label{eq:modifiedBilinearBiLaplace}
\end{align}
as well as
\begin{align}
\mathbf{  \tilde M}_{ij} &=& \sum\limits_{k \in \I_c} \sum\limits_{q=1}^{K} w_{q} \left( \varPhi_i \cdot \varPhi_j \,\sqrt{\det G}\right) (\xi^q,k) \label{eq:modifiedMassMatrix}  
\end{align}
and
\begin{align}
\mathbf{\tilde B}_j &=& \sum\limits_{k \in \I_c} \sum\limits_{q=1}^{K} w_{q} \left( (f \circ X) \cdot \varPhi_j \, \sqrt{\det G} \right)(\xi^q,k). \label{eq:modifiedRHS}
\end{align}
Now, instead of solving \eqref{eq:exactLinSystem}, we solve the system $\tilde \S \cdot \mathbf{U} = \mathrm{\tilde B}$.
The associated, modified variational problem reads as follows: Find $\tilde u_h \in V^0_h$ such that 
\begin{align}
 \tilde a(\tilde u_h,v_h) = \tilde \ell_f(v_h), \quad \forall ~v_h \in V_h^0,  \label{eq:modifiedWeakForm}
\end{align}
where  $\tilde a( \varphi_i, \varphi_j ) =\mathbf{\tilde S}_{ij}$ and $\tilde \ell_f(\varphi_j) = \mathbf{\tilde B}_{j}$.
If $\tilde a(.,.)$ is uniformly $V_h$-elliptic ($\mathbf{\tilde S}$ positive definite) existence and uniqueness of the discrete solution $\tilde u_h$ 
is ensured (cf. \cite{BrSc02}, for instance). 

\subsection*{Gaussian Quadrature.}
For second-order elliptic problems, standard error estimates for finite element schemes with Gaussian quadrature \cite[Chapter 4.1]{Ci02} 
imply that the expected order of convergence, in case of exact integration, is preserved if the quadrature scheme is exact for polynomials of degree $p=6$.
Transferring this to fourth-order problems we request exactness for polynomials of degree $p=4$.
This suggests to choose a Gaussian quadrature rule GA($p$) which guarantees this required exactness. 
For the Laplace-Beltrami equation \eqref{eq::strongLaplaceProblem} $K=12$ and for the surface bi-Laplacian equation \eqref{eq::strongBiLaplacianProblem} 
$K=6$ quadrature points have to be taking into account on the reference triangle $\cell$ 
(for symmetric Gaussian quadrature points on triangles see \cite{Du85}).

\subsection*{Adaptive Gaussian Quadrature.}
Special care is required close to EVs \cite{PeRe08}, because the basis functions $\Phi_i$ are no longer polynomials
and the second order derivatives are singular at the EV.  
Furthermore, the natural parametrization \cite{St99} produces only $C^0$-surfaces at EVs instead of $C^1$-parametrizations. 
To obtain a $C^1$-parametrization a suitable embedding has to be considered, e.g., using the \emph{characteristic map} \cite{BoZo04, PeRe08}.
For the implementation, the evaluation algorithm of \cite{St99} can still be used, thanks to the integral transformation rule.

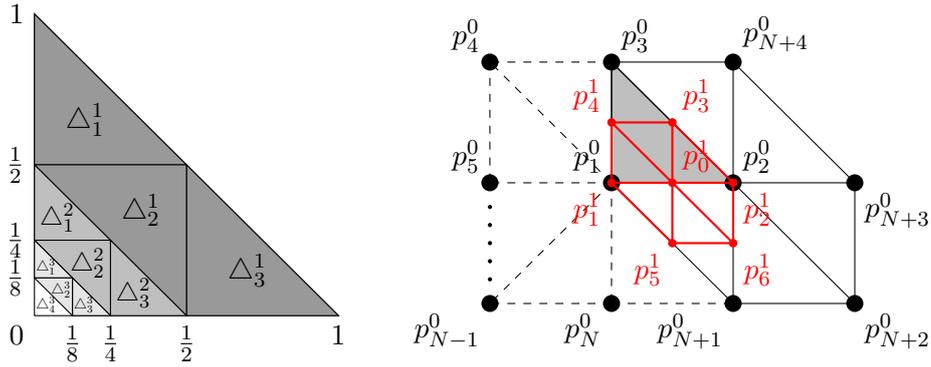
\begin{figure}
\begin{center}
\begin{tabular}{m{0.25\linewidth}m{0.5\linewidth}}
\begin{tikzpicture}[scale=4]
\draw[fill=white] (0,0) -- (0,1) -- (1,0) -- (0,0);
\draw[fill=light-gray1] (0,1/2) -- (0,1) -- (1/2,1/2) -- (0,1/2);
\draw[fill=light-gray1] (1/2,0) -- (1/2,1/2) -- (1,0) -- (1/2,0);
\draw[fill=light-gray1] (0,1/2) -- (1/2,0) -- (1/2,1/2) -- (0,1/2);
\draw[fill=light-gray2] (0,1/4) -- (0,1/2) -- (1/4,1/4) -- (0,1/4);
\draw[fill=light-gray2] (1/4,0) -- (1/4,1/4) -- (1/2,0) -- (1/4,0);
\draw[fill=light-gray2] (0,1/4) -- (1/4,0) -- (1/4,1/4) -- (0,1/4);
\draw (1/6,13/20) node {$\cell_1^1$};
\draw (7/20,7/20) node {$\cell_2^1$};
\draw (14/20,3/20) node {$\cell_3^1$};
\draw (3/40,13/40) node {\scalebox{.9}{$\cell_1^2$}};
\draw (7/40,7/40) node {\scalebox{.9}{$\cell_2^2$}};
\draw (13/40,3/40) node {\scalebox{.9}{$\cell_3^2$}};
%
\draw[fill=light-gray3] (0,1/8) -- (0,1/4) -- (1/8,1/8) -- (0,1/8);
\draw[fill=light-gray3] (1/8,0) -- (1/8,1/8) -- (1/4,0) -- (1/8,0);
\draw[fill=light-gray3] (0,1/8) -- (1/8,0) -- (1/8,1/8) -- (0,1/8);
\draw (3/80,13/80) node {\scalebox{.5}{$\cell_1^3$}};
\draw (7/80,7/80) node {\scalebox{.5}{$\cell_2^3$}};
\draw (13/80,3/80) node {\scalebox{.5}{$\cell_3^3$}};
\draw (3/80,3/80) node {\scalebox{.5}{$\cell_4^3$}};
%
\draw (0.125,0) node [below] {$\frac18$};
\draw (0.25,0) node [below] {$\frac14$};
\draw (0.5,0) node [below] {$\frac12$};
\draw (1,0) node [below] {$1$};

\draw (0,0.125) node [left] {$\frac18$};
\draw (0,0.25) node [left] {$\frac14$};
\draw (0,0.5) node [left] {$\frac12$};
\draw (0,1) node [left] {$1$};

\draw (0,0) node [below left] {$0$};
\end{tikzpicture} &
\quad \quad 
\begin{tikzpicture}[scale=1.6]
  \coordinate (x0) at (0,0); 
  \coordinate (x1) at (1,0); 
  \coordinate (x3) at (0,1); 
  \coordinate (x4) at (-1.,1); 
  \coordinate (x5) at (1.,-1); 
  \coordinate (x6) at (0.,-1); 
  \coordinate (x7) at (1.,1); 
  \coordinate (x8) at (2.,-1); 
  \coordinate (x9) at (2.,0.); 
  \coordinate (x10) at (.5,0.); %
  \coordinate (x11) at (-1.,0.); 
  \coordinate (x12) at (-1.,-1.); 
  \coordinate (x13) at (0.5,0.5); %
  \coordinate (x14) at (0,0.5); %
  \coordinate (x15) at (0.5,-0.5); %
  \coordinate (x16) at (1.,-0.5); %
  
  \draw[fill=light-gray2] (x0) -- (x1) -- (x3) -- (x0);
  \fill[black] (x3) circle (2.pt);
  \fill[black] (x4) circle (2.pt);
  \fill[black] (x5) circle (2.pt);
  \fill[black] (x6) circle (2.pt);
  \fill[black] (x7) circle (2.pt);
  \fill[black] (x8) circle (2.pt);
  \fill[black] (x9) circle (2.pt);
  \fill[black] (x11) circle (2.pt);
  \fill[black] (x12) circle (2.pt);
  \fill[black] (-1.,-0.20) circle (0.5pt);
  \fill[black] (-1.,-0.35) circle (0.5pt);
  \fill[black] (-1.,-0.5) circle (0.5pt);
  \fill[black] (-1.,-0.65) circle (0.5pt);
  \fill[black] (-1.,-0.80) circle (0.5pt);

  \draw (x0) node [above left] {$p^{0}_1$} -- (x3) node [above right] {$p^{0}_{3}$};
  \draw[dashed] (x0) -- (x4) node [above left] {$p^{0}_{4}$};
  \draw (x0) -- (x5) node [below left] {$p^{0}_{N+1}$};
  \draw[dashed] (x0) -- (x6) node [below left] {$p^{0}_{N}$};
  \draw (x1) node [above right] {$p^{0}_2$} -- (x3);
  \draw[dashed] (x3) -- (x4);
  \draw[dashed] (x6) -- (x5);
  \draw (x5) -- (x1);
  \draw (x3) -- (x7) node [above right] {$p^{0}_{N+4}$};
  \draw (x5) -- (x8) node [below right] {$p^{0}_{N+2}$};
  \draw (x1) -- (x7);
  \draw (x1) -- (x8);  
  \draw (x9) -- (x7);
  \draw (x9) -- (x8);
  \draw[dashed] (x4) -- (x11);
  \draw[dashed] (x6) -- (x12);
  \draw (x1) -- (x9) node [below right] {$p^{0}_{N+3}$};
  \draw[dashed] (x0) -- (x11) node [above left] {$p^{0}_{5}$};
  \draw[dashed] (x0) -- (x12) node [below left] {$p^{0}_{N-1}$};
  \fill (x0) circle (2pt);
  \fill (x1) circle (2pt);
  \fill[red] (x0) circle (1.pt);
  \fill[red] (x1) circle (1.pt);    
  \draw[red,thick] (x0) -- (x1);
  \draw[red,thick] (x0) node [below left] {\color{red}$p^{1}_{1}$} -- (x14) node [above left] {\color{red}$p^{1}_{4}$};
  \draw[red,thick] (x1) node [below right] {\color{red}$p^{1}_{2}$} -- (x13) node [above right] {\color{red}$p^{1}_{3}$};
  \draw[red,thick] (x13) -- (x14);
  \draw[red,thick] (x0) -- (x15) node [below left] {\color{red}$p^{1}_{5}$};
  \draw[red,thick] (x1) -- (x16) node [below right] {\color{red}$p^{1}_{6}$};
  \draw[red,thick] (x16) -- (x15);
  \draw[red,thick] (x10) node [above right] {\color{red}$p^{1}_{0}$} -- (x13);
  \draw[red,thick] (x10) -- (x14);
  \draw[red,thick] (x10) -- (x15);
  \draw[red,thick] (x10) -- (x16);
  \fill[red] (x10) circle (1.pt);
  \fill[red] (x13) circle (1.pt);
  \fill[red] (x14) circle (1.pt);
  \fill[red] (x15) circle (1.pt);
  \fill[red] (x16) circle (1.pt);
\end{tikzpicture}
\end{tabular}
\end{center}
\caption{Left: decomposition of the unit triangle $\cell$ in smaller triangles corresponding to the splines rings of the adaptive Gaussian quadrature AG($p,L$) for $L = 3$.
Right: illustration of the vertex configuration, with numbering as in Table \ref{tab:generalLookUp} (black straight and dashed lines) with a vertex $p_1^0$ of valence $N$  in the neighborhood of the edge $e$ connecting $p_1^0$ and $p_2^0$.
The red lines indicate the new vertices and edges in the 1-neighborhood of
$p_0^1$ after one level of subdivision.}
\label{fig:edgeMasks}
\end{figure}

Because of the structure of the subdivision scheme, the subdivision surface and correspondingly the basis functions $\Phi_i$ are spline functions on local triangular mesh rings 
around each EV \cite{PeRe08}. Thus, this subdivision ring structure can be used via an adaptive refinement strategy of the reference triangle $\triangle$ 
around the EV and an application of Gaussian quadrature on the resulting adaptive reference mesh to 
overcome the limitation of the standard Gaussian quadrature on the reference triangle.
More explicitly, for triangles with an EV we decompose the associated reference triangle $\cell$ into finer triangles $\cell_{i}^{l}$ ($l=1,\ldots, L$) of level $L$ (see Figure \ref{fig:edgeMasks}, left),
perform the corresponding Gaussian quadrature GA($p$) on all finer triangles $\cell_{i}^{l}$ and call this adaptive Gaussian quadrature AG($p,L$). 
Hence, the number of quadrature points $K$ on these adaptively refined triangles depends on  $L$: 
$K = (3 \cdot L + 1 ) \cdot 12$ for the Laplace-Beltrami equation \eqref{eq::strongLaplaceProblem} and $K = (3 \cdot L + 1 ) \cdot 6$ for the surface bi-Laplacian equation \eqref{eq::strongBiLaplacianProblem}. 
It is worth mentioning that this expensive scheme has to be applied only for a small number of triangles around the finitely many EVs.
Finally, let us remark that standard Romberg extrapolation does not lead to any improvement of the adaptive Gaussian quadrature due to the 
lack of smoothness of the integrands in the mesh size parameter and an adaption of the Romberg method for singular problems failed because the type of singularity is not explicitly known.

\subsection*{Barycenter Quadrature.}
Implementation of the adaptive Gaussian quadrature is a tedious issue. 
In particular, in graphics where the achievable maximal order of consistency is not needed, already the computing cost of the usual Gaussian quadrature is significant.
Therefore, reduced quadrature assembly is a common practice when implementing modeling or simulation tools based on 
NURBS as well as for subdivision surfaces (e.g. \cite{HuReSa10, CiOrSc00}).
The simplest and for subdivision surfaces widespread quadrature  is the barycentric quadrature (BC) with the center point $\xi^1 = (\frac13,\frac13)$ of the triangle $\cell$ 
and the weight $w_1 = \frac12$. 
This rule is applied to regular as well as to irregular triangles and integrates only affine functions exactly.
The method is also used in \cite{CiOrSc00} and leads to a reasonable consistency at least in the energy norm.


\subsection*{Mid-edge Quadrature.}
An alternative, which shows superior performance with respect to the achievable convergence rates (cf. Section \ref{sec5}), 
is the mid-edge quadrature (ME) with $K = 3$ quadrature points at the 
midpoints of the edges, i.e., $\xi^1 = (\frac12,0)$, $\xi^2 = (\frac12,\frac12)$ and $\xi^3 = (0,\frac12)$ with weights $w_1 = w_2 = w_3 = \frac16$.
The ME quadrature integrates exactly polynomials of degree $p=2$. 
Instead of following the direct implementation of \eqref{eq:modifiedBilinearLaplace}, \eqref{eq:modifiedBilinearBiLaplace}, \eqref{eq:modifiedMassMatrix}, and \eqref{eq:modifiedRHS}
we derive geometry-independent lookup tables of the basis function values and their derivatives on the unit triangle $\cell$ at midpoints of the edges.
This substantially simplifies the implementation.

Figure \ref{fig:edgeMasks} (left) depicts a sketch of a generic local control mesh with two interior vertices, one possibly EV $p^0_1$ of valence $N$ and one 
regular vertex $p^0_2$ of valence $6$. These two control points are surrounded by a fan of additional control points $p_i^0$ ($i=3,...,N+4)$. 
To perform one step of subdivision, we introduce new vertices (which are always regular) on every edge and update the positions of the current vertices. 
The new vertex positions are linear combinations of the old vertices $p_i^0$. Using the well-known position limit weights $\omega_i$ for regular vertices 
\cite{KoDaSe98} we can directly compute the limit position $ p^\infty_0 = \sum_{i=0}^{6} \omega_i p_i^1$ of the control vertex $p^1_0$ in dependence of the coarse mesh vertices $p_i^0$. 
This process applies to both the subdivision function values and their derivatives, which
are listed in Table \ref{tab:generalLookUp}.
More explicitly, for the hat function $\Lambda_i$ associated to a control vertex $p_i^0$, the value of the basis function 
$\varPhi_i(\xi^e,k)$ at the midpoint of the edge $e \in \I_e$ (connecting $p_1^0$ and $p_2^0$ with local coordinates $\xi^e$ in the triangle $(\cell,k)$ corresponding to the vertices
$p_1^0$,  $p_2^0$, and $p_3^0$) is a constant solely depending on $N$. 
The same holds for the derivatives $\varPhi_{i,1}$, $\varPhi_{i,2}$, $\varPhi_{i,11}$, $\varPhi_{i,12}$ and $\varPhi_{i,22}$ at $\xi^e$ in the directions on the reference triangle $\cell$.
Let us remark that, if $p_1^0$ and $p_2^0$ are EVs it is not clear if the global set of basis functions are linearly independent \cite{PeWu06a}.

Based on Table \ref{tab:generalLookUp}, the \emph{isogeometric subdivision approach} can be implemented by iterating over all edges retrieving values from these lookup tables without implementing 
the box spline basis functions, their derivatives and the complex subdivision process itself. The iteration over edges instead of elements avoids to evaluate basis function values twice.
Furthermore, the involved local matrices in the assembly of the mass and stiffness matrices are smaller for the mid-edge rule than for the barycenter rule.
More explicitly, in the regular case, the local IgA--matrices have $144$ entries for the BC rule versus $100$ for the ME rule.
Thus, based on the lookup tables, this leads to an overall faster assembly of mass and stiffness matrices than for the BC rule, even though the number of edges is $\tfrac32$ times the number of triangles on closed surfaces (e.g. see Table \ref{tab:assemblyTimes}).
Finally, let us remark that the mid-edge assembly process based on lookup tables can easily be generalized to other subdivision schemes like the Catmull--Clark \cite{CaCl78} or the Doo--Sabin scheme \cite{DoSa78}.\\

\renewcommand{\arraystretch}{1.5}
\begin{table}[th]
\footnotesize
\begin{center}
\begin{tabular}{|c|c|c|c|c|c|c|}
\hline
             &            $\varPhi_i(\xi^e,e)$                 &            $\varPhi_{i,1}(\xi^e,e)$       &             $\varPhi_{i,2}(\xi^e,e)$        & $\varPhi_{i,11}(\xi^e,e)$                & $\varPhi_{i,12}(\xi^e,e)$                             & $\varPhi_{i,22}(\xi^e,e)$  \\ \hline
 $i=1$     & $\frac{69 - 16 \cdot N \cdot \beta(N)}{192}$  & $\frac{-19 + 16 \cdot N \cdot \beta(N)}{24}$ & $\frac{-19 + 16 \cdot N \cdot \beta(N)}{48}$ & $\frac{5 - 16 \cdot N \cdot \beta(N)}{4}$ & $\frac{5-16 \cdot N \cdot \beta(N)}{8}$  &       $-1$       \\\hline
 $i=2$     & $\frac{62 + 16 \cdot \beta(N)}{192}$          & $\frac{14 - 16 \cdot \beta(N)}{24}$          & $\frac{14 - 16 \cdot \beta(N)}{48}$          & $\frac{-2 + 16 \cdot \beta(N)}{4}$        & $\frac{-1+8 \cdot \beta(N)}{4}$          &       $-1$       \\\hline
 $i=3$     & $\frac{25 + 16 \cdot \beta(N)}{192}$          & $\frac{1 - 16 \cdot \beta(N)}{24}$           & $\frac{19 - 16 \cdot \beta(N)}{48}$          & $\frac{-3 + 16 \cdot \beta(N)}{4}$        & $\frac{-3+16 \cdot \beta(N)}{8}$         &       $\frac12$ \\\hline
 $i=4$     & $\frac{2 + 16 \cdot \beta(N)}{192}$           & $\frac{-1 - 16 \cdot \beta(N)}{24}$          & $\frac{2 - 16 \cdot \beta(N)}{48}$           & $4 \cdot \beta(N)$                        & $\frac{-1+8 \cdot \beta(N)}{4}$          &       $0$        \\\hline
 $i\!=\!5, \cdots, N\!-\!1$     & $\frac{16 \cdot \beta(N)}{192}$               & $\frac{-2 \cdot \beta(N)}{3}$                & $\frac{-\beta(N)}{3}$                        & $4 \cdot \beta(N)$                        & $2\cdot \beta(N)$            &       $0$        \\\hline 
 $i={N}$   & $\frac{2 + 16 \cdot \beta(N)}{192}$           & $\frac{-1-16 \cdot \beta(N)}{24}$            & $\frac{-4 - 16 \cdot \beta(N)}{48}$          &       $4 \cdot \beta(N)$                  & $\frac{1+8 \cdot \beta(N)}{4}$          &       $\frac12$ \\\hline
 $i={N+1}$ & $\frac{25 + 16 \cdot \beta(N)}{192}$          & $\frac{1-16 \cdot \beta(N)}{24}$             & $\frac{-17 - 16 \cdot \beta(N)}{48}$         &  $\frac{(-3 + 16 \cdot \beta(N)}{4}$      & $\frac{-3+16 \cdot \beta(N)}{8}$          &       $\frac12$ \\\hline
 $i={N+2}$ & $\frac{3}{192}$                               & $\frac{1}{12}$                               & $\frac{-1}{48}$                              &       $\frac14$                           & $\frac{-1}{8}$                              &       $0$        \\\hline
 $i={N+3}$ & $\frac{1}{192}$                               & $\frac{1}{24}$                               & $\frac{1}{48}$                               &       $\frac14$                           & $\frac{1}{8}$                              &       $0$        \\\hline
 $i={N+4}$ & $\frac{3}{192}$                               & $\frac{1}{12}$                               & $\frac{5}{48}$                               &       $\frac14$                           & $\frac{3}{8}$                              &       $\frac12$ \\\hline
\hline \hline
 modification     &&&&&&  \\[-2ex]
   for $N=3$    &            $\varPhi_i(\xi^e,e)$                 &            $\varPhi_{i,1}(\xi^e,e)$       &             $\varPhi_{i,2}(\xi^e,e)$        & $\varPhi_{i,11}(\xi^e,e)$                & $\varPhi_{i,12}(\xi^e,e)$                             & $\varPhi_{i,22}(\xi^e,e)$  \\ \hline
 $i=3$     & $\frac{27 + 16 \cdot \beta(N)}{192}$ & $\frac{-2\beta(N)}{3}$ & $\frac{15  - 16 \cdot \beta(N)}{48}$ & $\frac{-3 + 16 \cdot \beta(N)}{4}$ & $\frac{-1+16 \cdot \beta(N)}{8}$ &   $1$ \\\hline
 $i=4$     & $\frac{27 + 16 \cdot \beta(N)}{192}$ & $\frac{-2\beta(N)}{3}$ & $\frac{-15 - 16 \cdot \beta(N)}{48}$ & $\frac{-3 + 16 \cdot \beta(N)}{4}$ & $\frac{-5+16 \cdot \beta(N)}{8}$ &    $\frac12$        \\\hline
\hline\hline
 modification     &&&&&&  \\[-2ex] 
 for $N=4$   &            $\varPhi_i(\xi^e,e)$                 &            $\varPhi_{i,1}(\xi^e,e)$       &             $\varPhi_{i,2}(\xi^e,e)$        & $\varPhi_{i,11}(\xi^e,e)$                & $\varPhi_{i,12}(\xi^e,e)$                             & $\varPhi_{i,22}(\xi^e,e)$  \\ \hline
 $i=3$     & $\frac{25 + 16 \cdot \beta(N)}{192}$ & $\frac{1 - 16 \cdot \beta(N)}{24}$  & $\frac{19 - 16 \cdot \beta(N)}{48}$  & $\frac{-3 + 16 \cdot \beta(N)}{4}$ & $\frac{-3+16 \cdot \beta(N)}{8}$ & $\frac12$ \\\hline
 $i=4$     & $\frac{4 + 16 \cdot \beta(N)}{192}$  & $\frac{-2 - 16 \cdot \beta(N)}{24}$ & $\frac{-2 - 16 \cdot \beta(N)}{48}$  & $4 \cdot \beta(N)$                    & $2 \cdot \beta(N)$                  & $\frac12$ \\\hline
 $i=5$     & $\frac{25 + 16 \cdot \beta(N)}{192}$ & $\frac{1 - 16 \cdot \beta(N)}{24}$  & $\frac{-17 - 16 \cdot \beta(N)}{48}$ & $\frac{-3 + 16 \cdot \beta(N)}{4}$ & $\frac{-3+16 \cdot \beta(N)}{8}$ & $\frac12$ \\\hline
    
\end{tabular}
\end{center}
\caption{Lookup tables for mid-edge quadrature rule ME. $N$ denotes the valence of vertex $p_1^0$ and $\beta(N) = \frac1N \left( \frac58 - \left( \frac38 + \frac28 \cos\left( \frac{2\pi}{N}\right) \right)^2 \right)$. 
For $N=3$ and $N=4$ and some indices $i$ the values differ from the those for general $N$. For the corresponding reference configuration see Figure \ref{fig:edgeMasks} (here $i$ corresponds to $p_i^0$).}
\label{tab:generalLookUp}
\end{table}

\renewcommand{\arraystretch}{1}

\begin{figure}[!ht]
\begin{center}
\includegraphics[scale=.25]{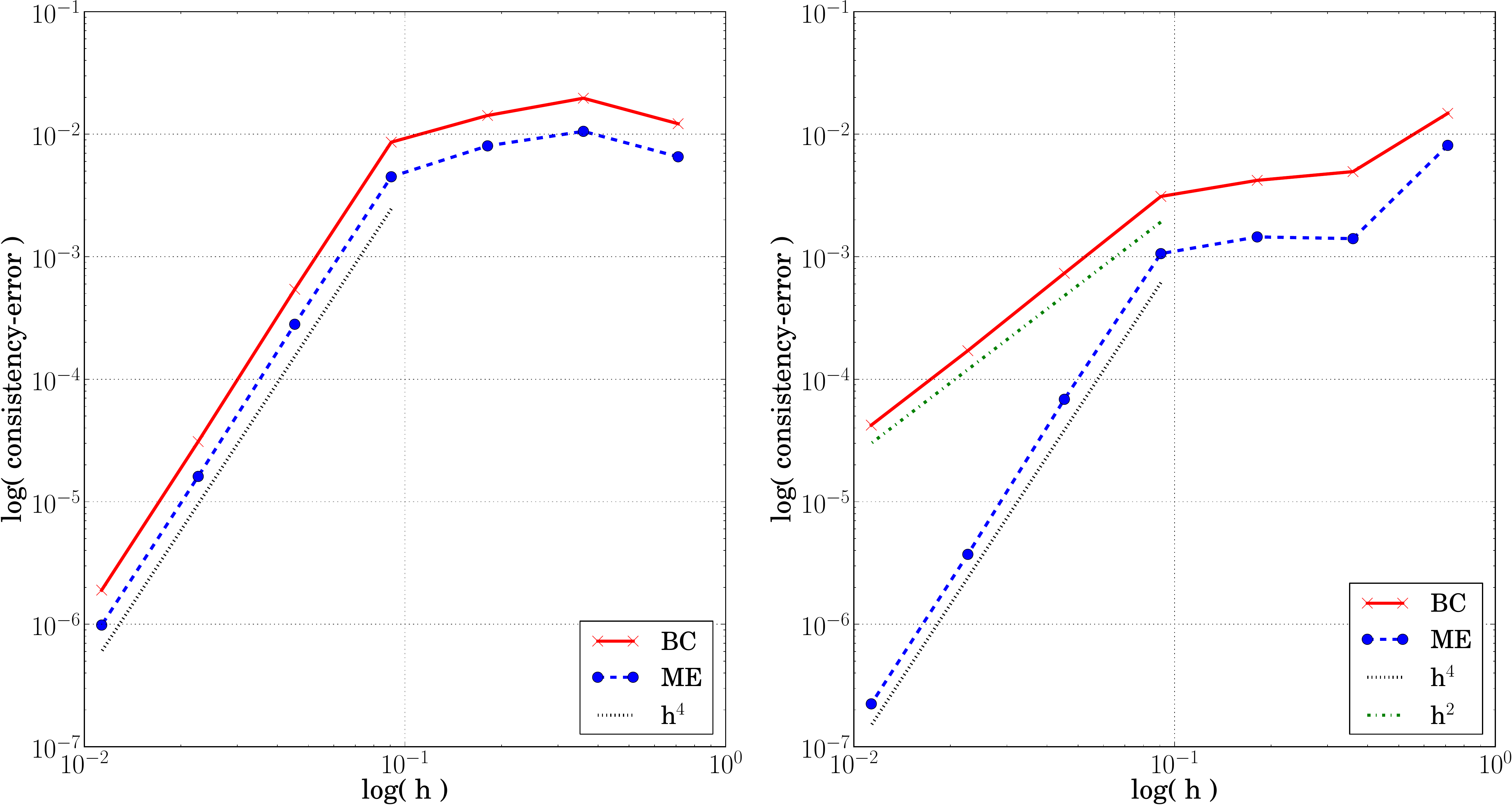}
\end{center}
\caption{The consistency error is shown in a log-log plot for the Laplacian problem on the torus (cf. Fig. \ref{fig:torus}) for varying grid size of the control mesh and for the 
barycenter quadrature (BC) and the midedge quadrature (ME).
}
\label{fig:torus_consistency}
\end{figure}

\noindent {\bf Remark}
Strang's Lemma (see \cite{StFi73,BrSc02}) provides the following error estimate for the numerical solution $\tilde u_h$ of \eqref{eq:modifiedWeakForm}
\begin{align*}
\| u - \tilde u_h \| &\leq \underbrace{\| u - u_h \|}_{\text{Approximation error (Cea's Lemma)}} + \underbrace{\| u_h - \tilde u_h \|}_{\text{Consistency error (Strang's Lemma)}} \\
&\leq \inf_{v_h \in V_h^0} \left( \| u - v_h \| + \sup_{w_h \in V_h^0} \frac{|a(v_h,w_h) - \tilde a(v_h,w_h)  |}{\|w_h\|} \right) + \sup_{w_h \in V_h^0} \frac{|\ell(w_h) - \tilde \ell(w_h)  |}{\|w_h\|}
\end{align*}
where $u$ is the continuous solution of \eqref{eq:generelWeakForm} and $u_h$ is the discrete solution of \eqref{eq:generelDiscreteWeakForm}.
The last two terms measure the \emph{consistency} of $\tilde a$ and $\tilde \ell_f$.
Here, $\| . \| := \| . \|_{H^1}$ for the Laplace-Beltrami problem \eqref{eq::strongLaplaceProblem} and $\| . \| := \| . \|_{H^2}$ for the surface bi-Laplacian problem \eqref{eq::strongBiLaplacianProblem}.
In fact, if a scheme with exact integration fulfills 
$\inf_{v_h \in V_h} \| u - v_h \| \leq C h^p$  (with $p \geq 1$ and a constant $C$),
we ask for a numerical quadrature that preserves this order, i.e.,
$$\sup_{w_h \in V_h} \left( \frac{|a(u_h,w_h) - \tilde a(u_h,w_h)  |}{\|w_h\|} + \frac{|\ell(w_h) - \tilde \ell(w_h)  |}{\|w_h\|} \right) \leq C h^p.$$
The optimal  $w_h = \sum_{i\in \I_v} w_i \varphi_i \in V_h^0$ for fixed $h$ is given by
$w_h = \tfrac{z_h}{\sqrt{a(z_h, z_h)}}$ where $z_h\in V_h^0$ is the solution of 
\begin{align*}
 a(z_h, \varphi_j) = a(u_h,\varphi_j) - \tilde a(u_h,\varphi_j), \quad \forall ~ j \in \I_v,
\end{align*}
with $\int_\M z_h \da = 0$ and $z_h = \sum_{i\in \I_v} z_i \varphi_i$. Figure \ref{fig:torus_consistency} 
plots the resulting consistency error depending on the mesh size $h$ of the control mesh. We observe an improved consistency 
by two orders for the mid-edge rule compared to the barycenter rule for the surface bi-Laplacian problem \eqref{eq::strongBiLaplacianProblem}.

\section{Numerical Applications} \label{sec5}
\begin{figure}[ht]
\begin{center}
\includegraphics[scale=1.]{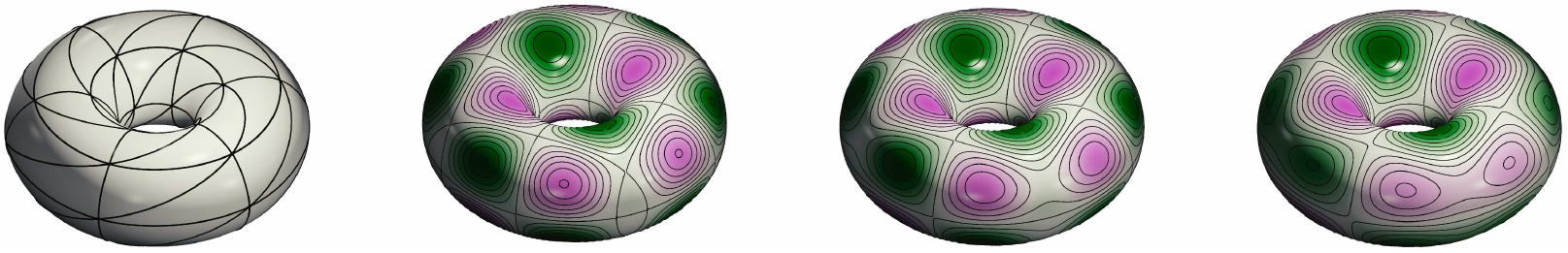} \\[2ex]
\includegraphics[scale=.32]{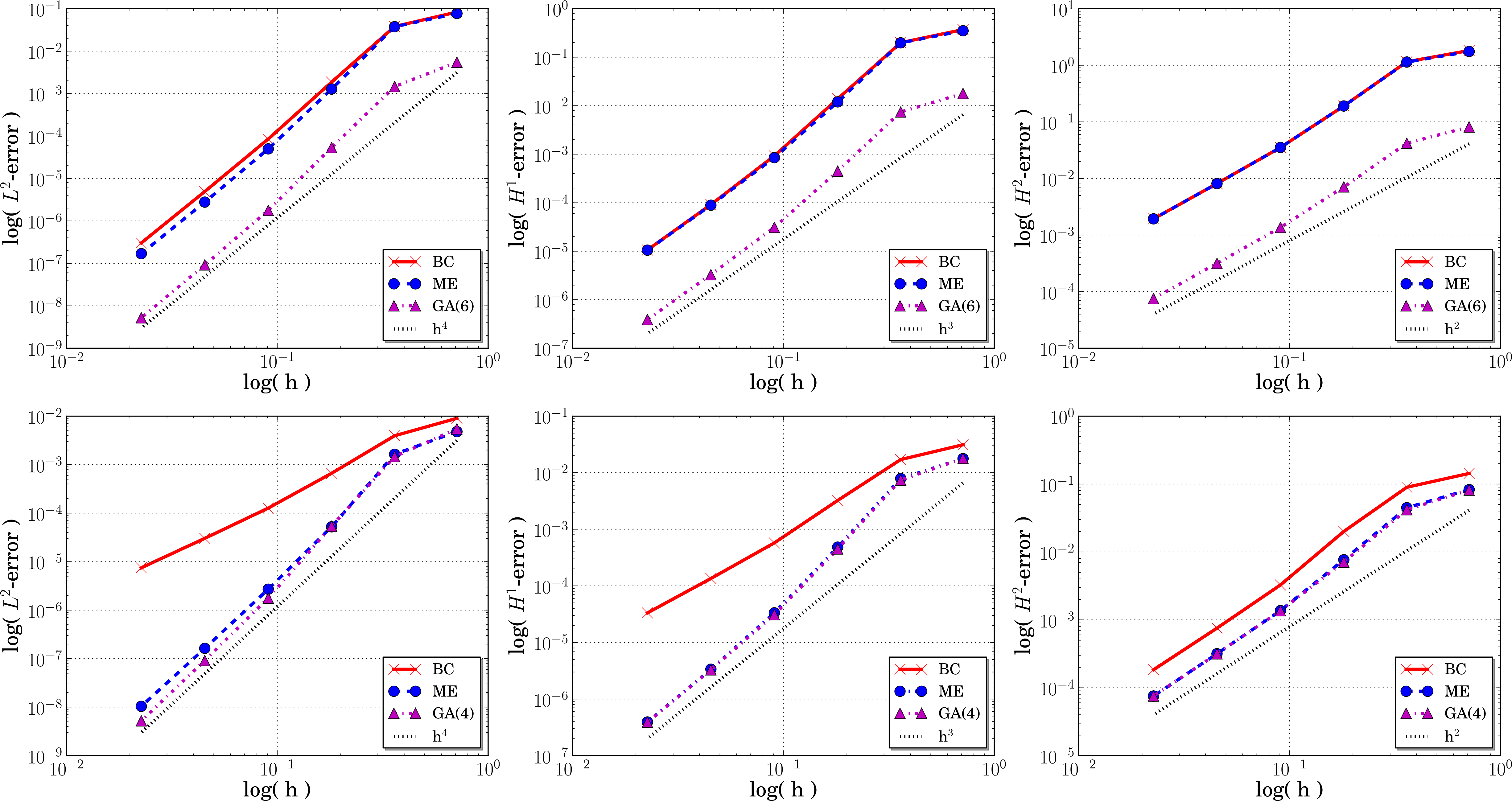}
\end{center}
\caption{Results for the torus: 
in the first row from left to right we depict the subdivision limit surface with control lines, isolines and color coding of the right hand side  
$f(y) = \sin(\pi y_1) \sin(\pi y_2) \sin(\pi y_3)$ ($-1.0$ \protect\resizebox{.06\linewidth}{!}{\protect\includegraphics[width=0.2\textwidth, height=20px]{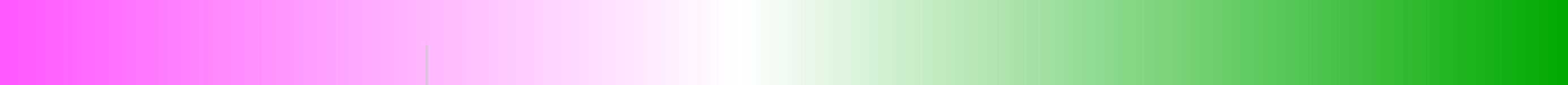}} $1.0$),
the (numerical reference) solution of the Laplace-Beltrami problem ($-0.054$ \protect\resizebox{.06\linewidth}{!}{\protect\includegraphics[width=0.2\textwidth, height=20px]{colorbar.pdf}} $0.054$),
and the (numerical reference) solution of the surface bi-Laplacian problem ($-0.003$ \protect\resizebox{.06\linewidth}{!}{\protect\includegraphics[width=0.2\textwidth, height=20px]{colorbar.pdf}} $0.003$) on $\M$;
in the second and third row log-log plots of the error are reported for the Laplace-Beltrami problem ($2^{nd}$ row) 
and the surface bi-Laplacian problem ($3^{rd}$ row), respectively (from left to right: $L^2$-norm, $H^1$-semi-norm and $H^2$-semi-norm).
}
\label{fig:torus}
\end{figure}
We have implemented the proposed methods in C++ and performed tests for four different subdivision surfaces: a torus surface with regular control mesh, 
a spherical surface (Spherical-3-4) with a control mesh with EVs of valence $3$ and $4$, 
a spherical surface (Spherical-5-12) with a control mesh with EVs of valence $5$ and $12$ and 
a complex real world hand model where the control mesh has altogether $119$ EVs of valence $4$, $5$, $7$, $8$, $9$ and $10$. 
Let us emphasize that, in the spirit of the isogeometric approach, the control mesh determining the limit subdivision geometry is kept fixed.
Furthermore, the limit surface differs from a torus (Figure \ref{fig:torus}) with circular centerline and cross section 
or a perfect sphere (Figure \ref{fig:Sphere12with5}). Then, we successively refine the control mesh using subdivision refinement to improve the accuracy of the discrete PDE solution.
To run the simulations for Spherical-3-4, Spherical-5-12 and the hand model the initial control mesh has to be subdivided at least once to avoid EVs in direct neighborhood.
In all tables, the mesh size $h$ refers to the mesh size of these control meshes.
\begin{figure}[!ht]
\begin{center}
\includegraphics[scale=1.]{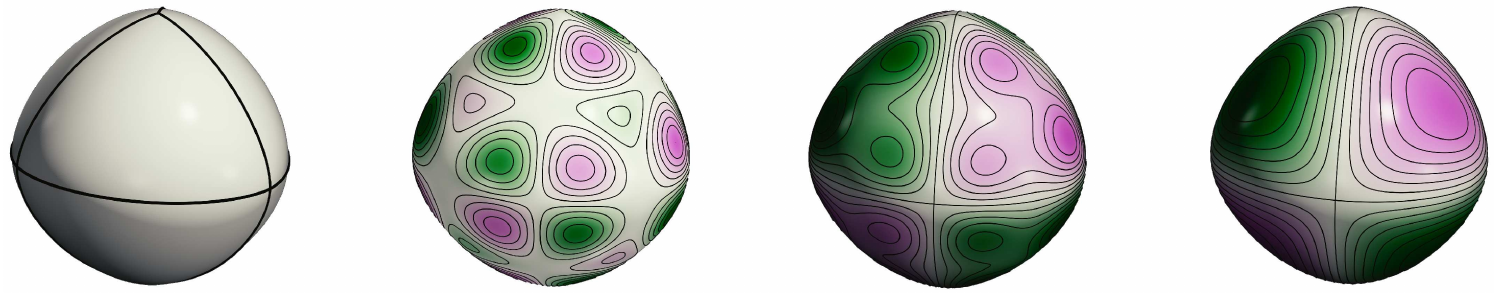} \\[2ex]
\includegraphics[scale=.32]{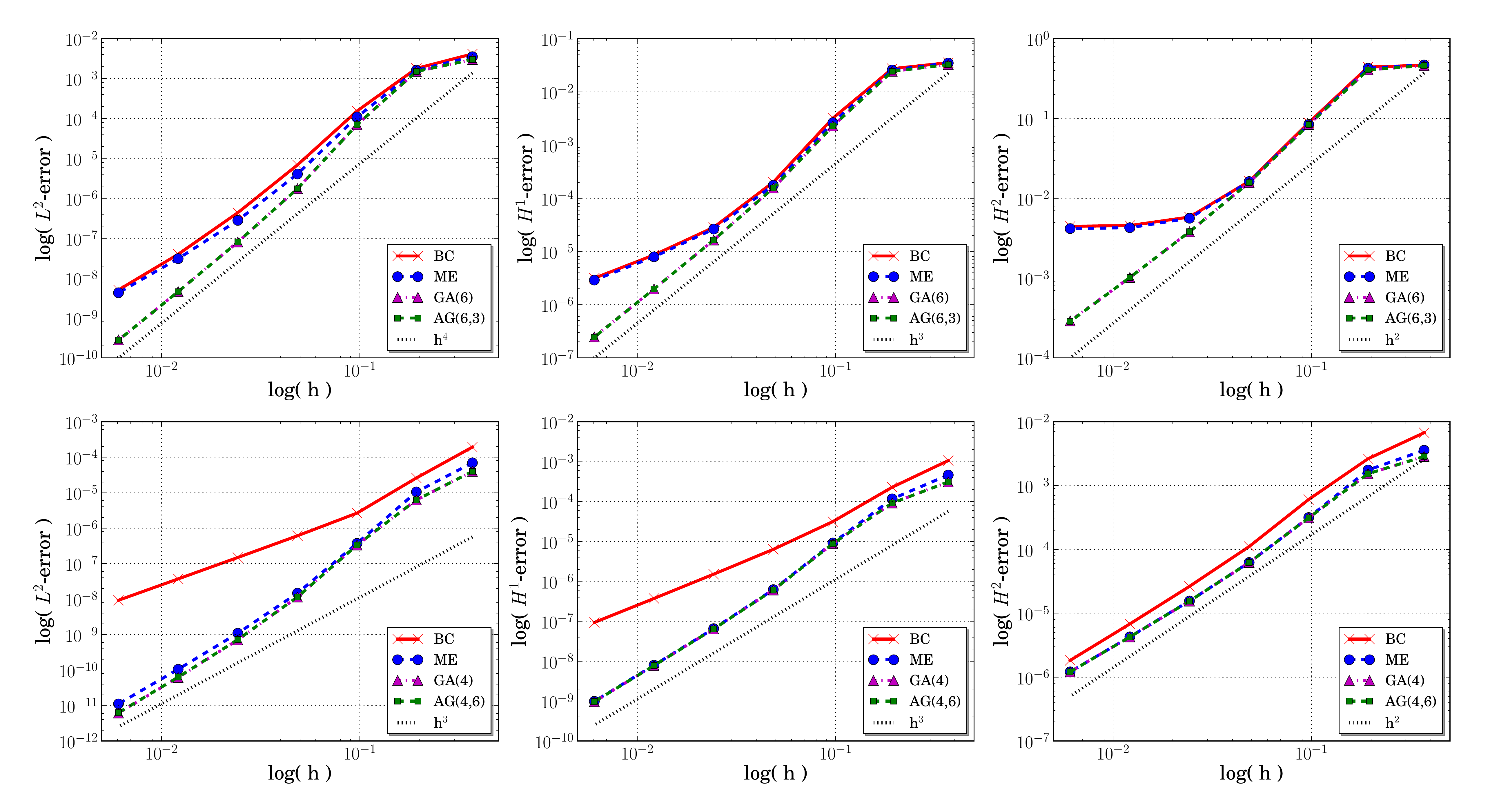}
\end{center}
\caption{Results for Spherical-3-4: 
As in Figure \ref{fig:torus} we show in the top row the limit surface, $f(y) = \sin(3 \pi y_1) \sin(3 \pi y_2) \sin(3 \pi y_3)$ ($-1.0$ \protect\resizebox{.06\linewidth}{!}{\protect\includegraphics[width=0.2\textwidth, height=20px]{colorbar.pdf}} $1.0$),
the solution of the Laplace-Beltrami problem ($-5.25e-3$ \protect\resizebox{.06\linewidth}{!}{\protect\includegraphics[width=0.2\textwidth, height=20px]{colorbar.pdf}} $5.25e-3$),
and the solution of the surface bi-Laplacian problem ($-9.78e-5$ \protect\resizebox{.06\linewidth}{!}{\protect\includegraphics[width=0.2\textwidth, height=20px]{colorbar.pdf}} $9.78e-5$),
together with the error-plots for the Laplace-Beltrami problem ($2^{nd}$ row) and the surface bi-Laplacian problem ($3^{rd}$ row) (from left to right: $L^2$-norm, $H^1$-semi-norm and $H^2$-semi-norm).
}
\label{fig:Sphere3with4}
\end{figure}
\begin{figure}[!ht]
\begin{center}
\includegraphics[scale=1.]{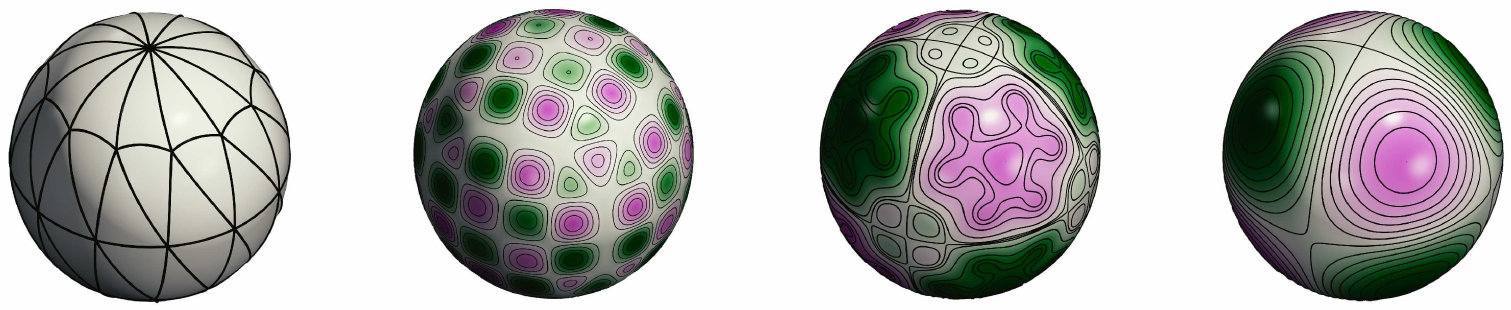} \\[2ex]
\includegraphics[scale=.32]{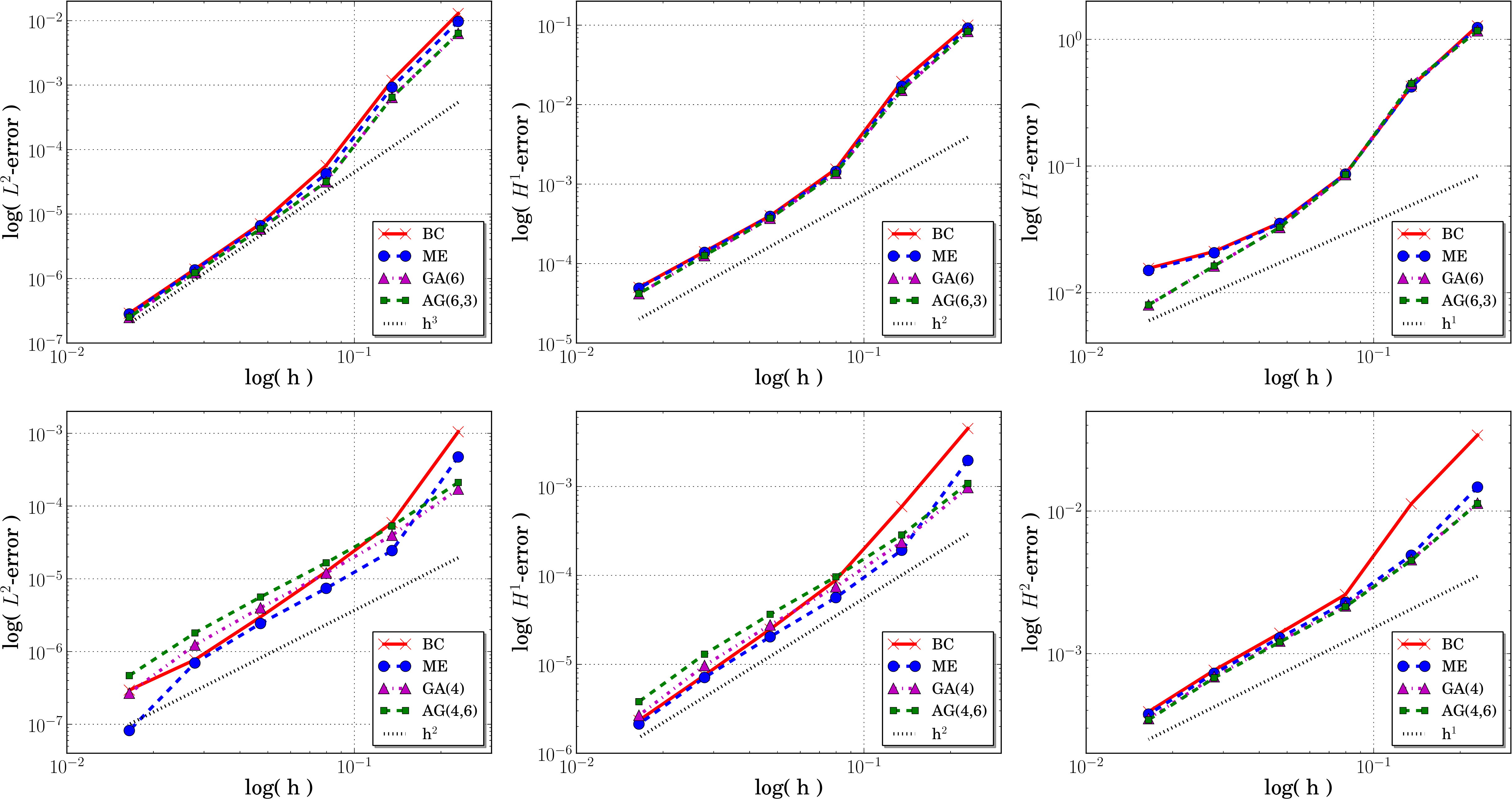}
\end{center}
\caption{Results for Spherical-5-12: 
We follow the presentation in Fig. \ref{fig:Sphere3with4} and plot the limit surface with IgA-element lines, the (same) right hand side function $f$ now on this surface, 
the solution of the Laplace-Beltrami problem ($-1.58e-2$ \protect\resizebox{.06\linewidth}{!}{\protect\includegraphics[width=0.2\textwidth, height=20px]{colorbar.pdf}} $1.58e-2$), and
the solution of the surface bi-Laplacian problem ($-1.27e-3$ \protect\resizebox{.06\linewidth}{!}{\protect\includegraphics[width=0.2\textwidth, height=20px]{colorbar.pdf}} $1.27e-3$),
again together with the error-plots for the Laplace-Beltrami problem ($2^{nd}$ row) and the surface bi-Laplacian problem ($3^{rd}$ row) (from left to right: $L^2$-norm, $H^1$-semi-norm and $H^2$-semi-norm).
}
\label{fig:Sphere12with5}
\end{figure}

On all of these surfaces we solve the two model problems \eqref{eq::strongLaplaceProblem} and \eqref{eq::strongBiLaplacianProblem} 
for a given right-hand side $f$ and compare the asymptotic error for different norms: the $L^2$-norm, the $H^1$-semi-norm and the $H^2$-semi-norm. 
To evaluate the \emph{experimental order of convergence} (eoc), we have computed a reference solution for a control mesh with one additional level of global refinement compared to the finest mesh.
Furthermore to compute the reference solution we used the adaptive Gaussian quadrature AG($p$,$L$) with $L=3$ for \eqref{eq::strongLaplaceProblem} and for $L=6$ for \eqref{eq::strongBiLaplacianProblem} 
(the estimated convergence rates did not change for larger $L$). As a consequence, the error plots for the finest discretization are less reliable as it becomes apparent in Figure \ref{fig:Sphere12with5} 
for the $L^2$-error plot for problem \eqref{eq::strongBiLaplacianProblem} and the spherical shape with low valence EVs (Spherical-3-4).
For all the experiments, double precision arithmetic is used.

Figure \ref{fig:torus} shows our results for the torus with regular control mesh without EVs.
All quadrature rules achieve optimal convergence rates (cf. \cite{Ko90}) for the second-order problem \eqref{eq::strongLaplaceProblem}, i.e., $4$ in the $L^2$-norm, 
$3$ in the $H^1$-semi-norm and $2$ in the $H^2$-semi-norm (with non-adaptive Gaussian quadrature GA($p$)). 
Here, ''optimal'' reflects what we expect for the quartic box spline \cite{Ko90}.
For the fourth-order problem \eqref{eq::strongBiLaplacianProblem} all quadrature rules achieve optimal rates 
 $4$ in the $L^2$-norm,  $3$ in the $H^1$-semi-norm and  $2$ in the $H^2$-semi-norm (with non-adaptive Gaussian quadrature GA($p$)), except the BC rule. 
Let us remark that the $C^2$-regularity allows for an Aubin-Nitzsche type argument \cite{Ci02} to obtain estimates in norms other than the energy norm. 
In all cases GA($p$) performed best.

The numerical results of the spherical control mesh with EVs of valence $3$ and $4$ (Spherical-3-4) depicted in Figure \ref{fig:Sphere3with4}. 
For the second-order PDE \eqref{eq::strongLaplaceProblem} Gaussian quadrature and adaptive Gaussian quadrature achieve the optimal order of convergence. 
On finer resolutions ME and BC do not achieve optimal convergence rates. GA, AG and ME achieve similar behavior for the fourth-order problem and again the BC method performs worse.
The reduced convergence rate in the $L^2$-norm for the fourth-order problem is expected to reflect the reduced regularity $C^1 \cap H^2$, which precludes the application of an Aubin-Nitzsche type argument.

For the control mesh with  EVs of valence $5$ and $12$ (Spherical-5-12) the convergence results are depicted in 
Figure \ref{fig:Sphere12with5} and coincide with our general findings for control meshes with EVs of valence greater than $6$.
Here, all quadrature schemes show a similar performance,
for the Laplace-Beltrami problem \eqref{eq::strongLaplaceProblem}:  $3$ in the $L^2$-norm, $2$ in the $H^1$-semi-norm and $1$ in the $H^2$-semi-norm, 
and for the surface bi-Laplacian problem \eqref{eq::strongBiLaplacianProblem} : $2$ in $L^2$-norm, $2$ in the $H^1$-semi-norm and $1$ in the $H^2$-semi-norm.
This observed loss of approximation order compared to the quartic box spline coincides with the theoretical work by \cite{Ar01} with the 
reasoning that Loop subdivision functions cannot reproduce cubic polynomials around extraordinary vertices of valence greater then $6$.

\begin{figure}[!ht]
\begin{center}
\includegraphics[scale=1.]{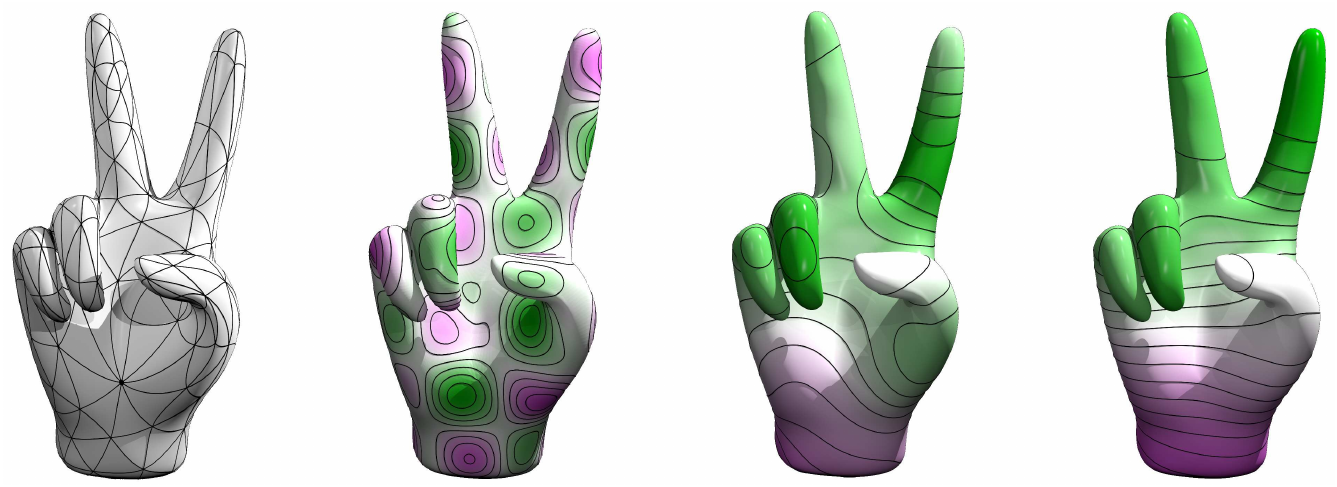} \\[2ex]
\includegraphics[scale=.32]{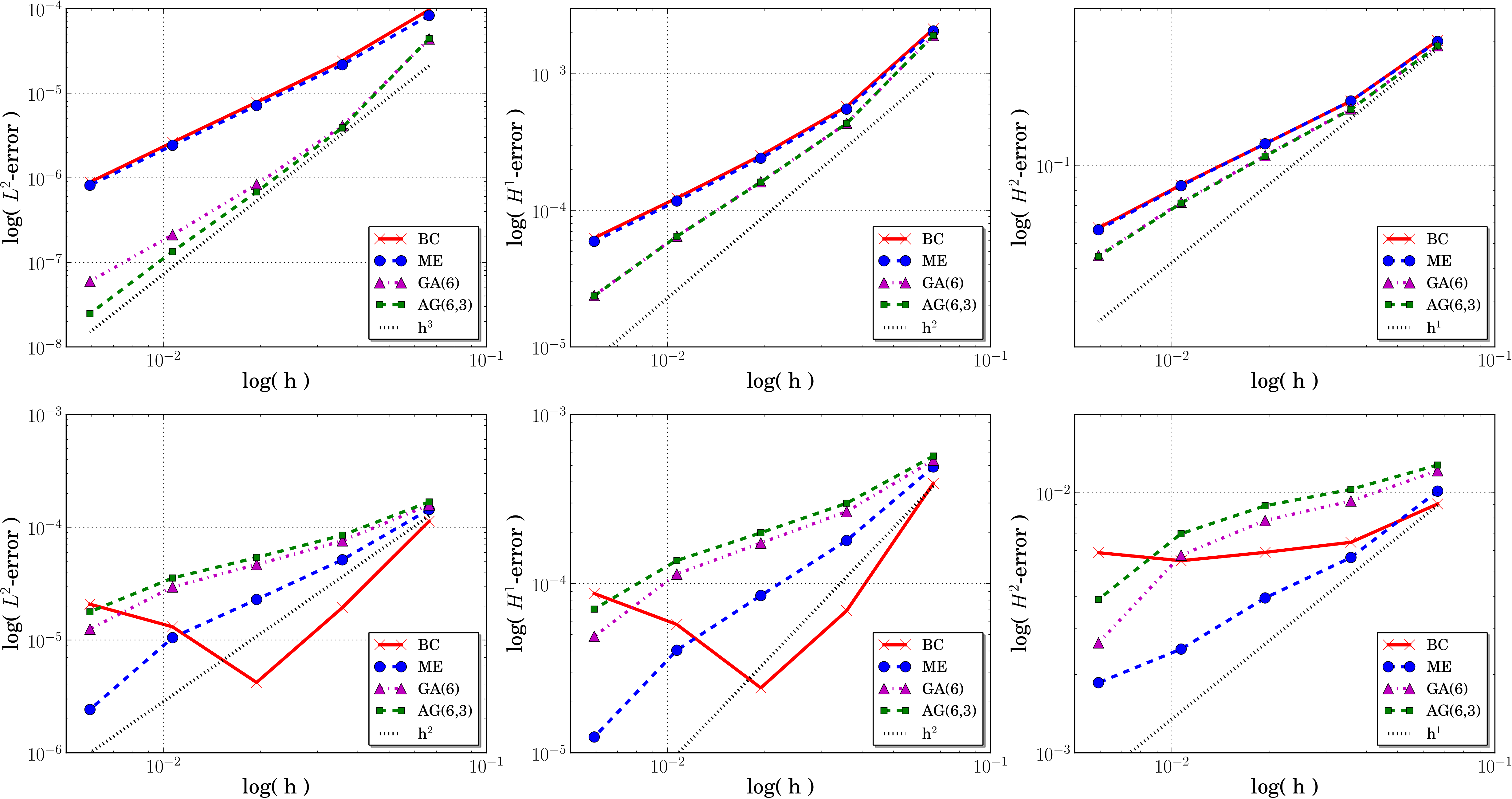}
\end{center}
\caption{Results for complex hand shape: 
We follow the presentation in Fig. \ref{fig:Sphere3with4} and plot the limit surface with IgA-element lines, the (same) right hand side function $f$ now on this surface, 
the solution of the Laplace-Beltrami problem ($-1.47e-2$ \protect\resizebox{.06\linewidth}{!}{\protect\includegraphics[width=0.2\textwidth, height=20px]{colorbar.pdf}} $1.14e-2$) and
the solution of the surface bi-Laplacian problem ($-8.14e.4$ \protect\resizebox{.06\linewidth}{!}{\protect\includegraphics[width=0.2\textwidth, height=20px]{colorbar.pdf}} $7.56e.4$),
again together with the error-plots for the Laplace-Beltrami problem ($2^{nd}$ row) and the surface bi-Laplacian problem ($3^{rd}$ row) (from left to right: $L^2$-norm, $H^1$-semi-norm and $H^2$-semi-norm).
}
\label{fig:handB}
\end{figure}

Finally, we consider the hand model as a complex subdivision surface with many EVs in Figure \ref{fig:handB}.
Given the complexity of the geometry and the number and valence of the extraordinary vertices both for the Laplace-Beltrami 
and surface bi-Laplacian problem the asymptotic regime of error reduction seems not to be reached even though the reference 
solution is computed on a mesh with $750 k$ vertices. 
The BC fails for Bi-Laplacian probably caused by the increased number of EVs and their partially high valences.
Figure \ref{fig:handB_eigen} shows eigenfunctions of the Laplace-Beltrami operator computed via inverse vector iteration with projection.
The depicted results underline that the methods discussed so far are
also applicable in the numerical eigenmode analysis, which turned out to be an indispensable tool in geometric data analysis and modeling.

\renewcommand{\arraystretch}{1.0}
\begin{table}
\begin{center}
\small
\begin{tabular}{|cc|cc|cc|cc|cc|}
\hline
\multicolumn{2}{|c}{} & \multicolumn{2}{|c|}{GA($p$)} & \multicolumn{2}{|c|}{AG($p$,$L$)} & \multicolumn{2}{|c|}{BC} & \multicolumn{2}{|c|}{ME} \\
Geometry    & $\mid\I_v\mid$ & $-\LBO$ & $(-\LBO)^2$ & $-\LBO$ & $(-\LBO)^2$ & $-\LBO$ & $(-\LBO)^2$ & $-\LBO$ & $(-\LBO)^2$  \\ \hline
Spherical-5-12 &   9218    & 1.209   & 1.417   & 1.374  & 1.563   & 0.164 &  0.272   & 0.107 & 0.136 \\
Hand        &   11586   & 1.914   & 2.233   & 2.615  & 3.424   & 0.269 &  0.460   & 0.157 & 0.198 \\
\hline
\end{tabular}
\end{center}
\caption{Comparison of assembly times for the stiffness matrices \eqref{eq:modifiedBilinearLaplace} and \eqref{eq:modifiedBilinearBiLaplace} 
and the different quadrature rules (Gaussian GA($p$), adaptive Gaussian AG($p$,$L$), barycenter (BC) and mid-edge (ME)).
All times reported in seconds.}
\label{tab:assemblyTimes}
\end{table}
\renewcommand{\arraystretch}{1.}
Compared to Gaussian quadrature and in particular to the adaptive Gaussian quadrature the mid-edge and the barycenter quadrature are significantly cheaper as reported in Table \ref{tab:assemblyTimes}. 
Because of our implementation, the mid-edge rule based on lookup tables and assembly via an edge-iterator performs even better than the (non-optimized) barycenter rule.

\section{Conclusions} \label{sec6}
We have explored the discretizations of geometric partial differential
equations on subdivision surfaces by means of an isogeometric approach.
To this end, we have focused on Loop's subdivision method and have investigated 
the impact of different numerical quadrature rules on the expected convergence rates.
As expected, there is a trade-off between robustness and computational effort. 
The adaptive Gaussian quadrature turns out to be the most robust in the carefully selected set of test cases
but at the same time computationally demanding for instance for real-time applications in modeling and animation.
In fact, the mid-edge quadrature  based on lookup tables has been singled out 
as a promising compromise for both graphics
and engineering applications. It is very easy to implement and
it performed well in all considered test cases.  
Fourth-order problems like the surface bi-Laplacian problem considered here, are by far more critical with respect to the accuracy of the numerical approximation 
than second-order problems.
In general, adaptivity around EVs is essential
to ensure the overall robustness on very fine computational meshes.


\section*{Acknowledgement}
We acknowledge support by the FWF in Austria under the grant S117 (NFN) and the DFG
in Germany under the grant Ru 567/14-1.

%
\bibliographystyle{elsarticle-num}
\bibliography{own,all,newBibTexItems}

\end{document}